\documentclass[12pt]{amsart}
\usepackage{amsfonts}
\usepackage{amssymb}
\usepackage{amscd}

\newcommand{\F}{{\mathbb F}_2}
\newcommand{\St}{{\bf S}_3}
\newcommand{\Dt}{{\bf D}_2}
\newcommand{\Af}{{\bf A}_4}
\newcommand{\Ps}{PSL_2({\mathbb Z})}

\newcommand{\calO}{{\mathcal O}}

\newcommand{\Z}{\mathbb{Z}}

\newcommand{\prf}{\medskip \noindent {\it Proof.} \quad}
\newcommand{\fin}{\quad $\square$}
\newcommand{\ararrow}{\buildrel \alpha \over \longrightarrow}
\newcommand{\brarrow}{\buildrel \rho \over \longrightarrow}
\newcommand{\drarrow}{\buildrel \delta \over \longrightarrow}

\newcommand{\G}[1]{{\Gamma}_{#1}}

\newtheorem{theorem}{Theorem}[section]
\newtheorem{lemma}[theorem]{Lemma}

\theoremstyle{definition}
\newtheorem{definition}[theorem]{Definition}

\theoremstyle{remark}
\newtheorem{remark}[theorem]{Remark}
\numberwithin{equation}{section}

\begin{document}

\title{The Integral Cohomology of the Bianchi Groups}

\author{Ethan Berkove}
\address{Lafayette College\\
 Department of Mathematics \\
 Easton, PA 18042}
\email{berkovee@lafayette.edu}


\begin{abstract}
We calculate the integral cohomology ring structure for various members of the Bianchi group family.  The main tools we use are the Bockstein spectral sequence and a long exact sequence derived from Bass-Serre theory. 
\end{abstract}

\maketitle

\section{Introduction}
In the calculation of $H^*(G;M)$, the cohomology of a group $G$ with coefficients in $M$, it is standard practice to let $M$ be a field in order to simplify the mathematics.  The integral cohomology ring, with $M = \Z$, can be much harder to determine, particularly when the group contains a lot of torsion.  However, it many cases it is still possible to get a complete answer.  In this paper we calculate the integral cohomology rings for a number of Bianchi groups.  These groups are particularly amenable to calculation as their torsion occurs only at the primes 2 and 3, their mod-2 cohomology rings and integral homology groups are known, and their 3-primary cohomology is easily calculable.  The groups themselves can be built in stages out of finite groups using amalgamated products and HNN extensions.  Therefore, their group cohomology is closely related to the cohomology of these finite pieces.   This paper completes a calculation that has partially appeared in a number of places:   For various members of the Bianchi group family, the integral homology (additive structure) can be found in \cite{S-V}, the rational homology in \cite{V}, and the mod-2 cohomology in \cite{Be}.

The main tools we use at the prime 2 are the Bockstein spectral sequence and long exact sequences that come from Bass-Serre theory.  At the prime 3 we perform the calculations solely using long exact sequences.  This paper is organized as follows.  In the next section we describe the Bianchi groups and the techniques we will need to determine their integral cohomology rings.  In the third section we apply these techniques to determine some intermediate results.  In the fourth section we give complete details for the calculation of $H^*(G;\Z)$ where $G = \G2, \G6$.  In the final section we 
do the cases $G = \G1, \G3, \G5, \G7, \G{10}, \G{11}$.

\section{Background}
For $d$ a square-free integer, let $\calO_d$ denote the ring of algebraic integers in the imaginary quadratic extension field of the rationals, ${\mathbb Q}(\sqrt{-d})$.  As a $\Z$-module, $\calO_d \cong \Z \oplus \Z \ \omega$, where $\omega = \sqrt{-d}$ when $d \equiv 3$ mod $4$ and  $\omega = \frac{1 + \sqrt{-d}}{2}$ when $d \equiv 1,2$ mod $4$.  The Bianchi groups $\Gamma_d$ are the projective special linear groups $PSL_2(\calO_d)$, and can be thought of as the natural generalizations of the modular group, $\Ps \cong \Z/2 \ast \Z/3$.  

The Bianchi groups are excellent candidates for integral cohomology calculations.  Torsion in Bianchi groups occurs at the primes 2 and 3, so one only needs to perform two $p$-primary calculations.  Furthermore, the most complicated torsion occurs at the prime $p=2$, and mod-2 cohomology rings (as modules over the Steenrod algebra) have already been calculated for the Bianchi groups $d = 1, 2, 3, 5, 6, 7, 10$ and $11$ in \cite{Be}. 

Starting with mod-$p$ cohomology information, the $p$-primary component of the integral cohomology can be reconstructed using a handy tool: the Bockstein spectral sequence.  This spectral sequence, which first appeared in \cite{Browd}, is associated to the following exact couple.  
$$ \ldots {\buildrel \delta \over \longrightarrow}  H^*(X;\Z) {\buildrel {(\times p^*)} \over \longrightarrow} H^*(X;\Z) {\buildrel red_{p^*} \over \longrightarrow} H^*(X;{\mathbb F}_p) {\buildrel \delta \over \longrightarrow}  H^{*+1}(X;\Z) {\buildrel {(\times p^*)} \over \longrightarrow}\ldots $$


The first differential is $d_1  = red_{p^*} \circ \delta$.  In fact, $d_1 = \beta$, the classical Bockstein homomorphism, which can be identified with the Steenrod square $Sq^1$ when $p=2$.  Therefore, if $Sq^1 (w) \neq 0$ for some class $w \in H^n(X;\F)$, then $\delta (w)$ represents some cohomology class  $\hat{w} \in H^{*+1}(X;\Z)$ of order 2.  Elements that appear in the integral cohomology at different pages of the spectral sequence have different orders; the $r^{th}$ differential $d_r$ can be identified with the $r^{th}$ order Bockstein homomorphism \cite{McC} 
$$\beta_r :H^*(X; \Z /2^r \Z) {\buildrel \times 2^{r-1} \over \longrightarrow} H^{*+1}(X; \Z /2^r \Z).$$
Thus, elements that survive to the $E_r$ page have order at least $2^r$.
The Bockstein spectral sequence converges to $(H^*(X;\Z) / \text{$p$-torsion} ) \otimes_{\Z} {\mathbb F}_p$.  We use this sequence to calculate group cohomology via the identification $H^*(G;M) = H^*(BG;M)$, where $BG$ is the classifying space of $G$.  

For the Bianchi groups we consider, the Bockstein spectral sequence applied to the calculations in \cite{Be} almost collapses at the $E_1$ page, with at most one torsion class surviving to the $E_2$ page.  Thus, almost all 2-torsion in the resulting integral cohomology rings has order 2, a result consistent with the calculations in  \cite{S-V}.  Although it would be convenient if this were the only tool needed in this paper, the Bockstein spectral sequence is less effective at identifying torsion-free classes, where knowledge of all Bockstein homomorphisms is required.  Further, we need a technique to determine the 3-primary cohomology ring.  

Fortunately, there is another technique.  The Bianchi groups we consider in this article are built in stages using amalgamated products and HNN extensions from copies of $\Z$ and five finite subgroups:  the cyclic groups of order two and three ($\Z/2$ and $\Z/3$), the Klein four group ($\Dt \cong \Z/2 \times \Z/2$), the symmetric group on three letters ($\St$), and the alternating group on four letters ($\Af$) \cite{G-S}.    We show how the cohomology of a group built via these constructions depends on the cohomology of its component subgroups.  This then allows us to determine the 3-primary and rational cohomology rings.  At the prime 2, this approach is also helpful for identifying the few classes of order 4 in the 2-primary cohomology ring.  We summarize some useful results, starting with the definition of an HNN extension.  A more thorough treatment can be found in \cite{Be}.
\begin{definition}Let $G_1$ be a group, $H$ be a subgroup and 
$\theta:H \rightarrow G_1$ be a monomorphism.  Then an {\it HNN extension 
of $G_1$} is a group $$\Gamma = \{ t,G_1 \mid t^{-1}at = \theta(a), a \in H \}.$$  
$\Gamma$ is denoted by $G_1 \ast_H$.  $G_1$ is called the {\it base} and 
$H$ is called the {\it associated subgroup}. 
\end{definition}

When $\Gamma \cong G_1 \ast_H G_2$ is an amalgamated product or $\Gamma \cong G_1 \ast_H$ is an HNN extension, Bass-Serre theory \cite{S} provides a long exact sequence.
\begin{equation}
\cdots \ararrow  H^{j-1}(H)  
       \drarrow H^j(\Gamma) \brarrow \oplus H^j(G_i)  
       \ararrow  H^j(H) \drarrow \cdots
\label{LES}
\end{equation}
The maps above can be explicitly described.  The map $\rho$ is a sum of restriction maps to each summand.  When $\Gamma$ is an amalgamated product, $\alpha$ is the difference of restriction maps, $\alpha = res^{G_1}_H - res^{G_2}_H$.  When $G$ is an HNN extension, the ``twisting'' induced by $\theta$ comes into play, and $\alpha = res^{G_1}_{H} - \theta^*$.  

\begin{remark}
Whenever $\Gamma$ an HNN extension, $\rho$ is an isomorphism in degree zero.  Consequently, $\alpha$ factors through zero and $\delta$ is an isomorphism, generating a torsion-free class in degree 1.  We call such classes ``HNN classes.''  By degree considerations, they are exterior.
\end{remark}

The map $\rho$ respects cup products, since it is a restriction map in the case of an HNN extension and a sum of restrictions in the case of an amalgamated product.  Thus, by exactness of $H^*(\Gamma) \brarrow \oplus H^*(G_i) \ararrow H^*(H)$, we have that $im(\rho) \cong ker(\alpha)$, and this is an {\bf isomorphism of rings}.  In most of our calculations the map $\alpha$ is surjective, so the long exact sequence of Formula \ref{LES} breaks into short exact sequences.  Hence the ring isomorphism determines most of the cup product structure of $H^*(\Gamma)$, which can be used to determine cup products in the cases where the Bockstein spectral sequence information is incomplete or ambiguous.

In those cases where $\alpha$ is not surjective, classes can also arise in $H^*(\Gamma)$ from the image of $\delta$.  The following lemmas describe how to determine products with these classes.

\begin{lemma} \cite{Sp}.  For $u \in H^p(\Gamma)$ and $v \in H^q(H),\ \delta (res^{\Gamma}_H u \cup v) = u \cup \delta (v),$ and $\delta (v \cup res^{\Gamma}_H u) = \delta (v) \cup u$.  \label{compat} 
\end{lemma}

\begin{lemma} Classes that arise from the image of $\delta$ have trivial products with each other. \label{PY}  
\end{lemma}

\prf See Theorem 4.8 in \cite{P-Y}  \fin

Finally, on occasion we will need to determine restriction maps on integral cohomology.  If the restriction map is known at the level of mod-2 cohomology, one can extend the map by appealing to the following long exact sequences of commuting squares.  
\hfill

\bigskip
\setlength{\unitlength}{1cm}
\begin{picture}(10,2.5)
  \put(2,2){$\ldots \longrightarrow H^*(G) {\buildrel {\times p} \over \longrightarrow} H^*(G) {\buildrel red_{p^*} 
      \over \longrightarrow} H^*(G;{\mathbb F}_p) {\buildrel \delta \over \longrightarrow} \ldots$}
  \put(2,0){$\ldots \longrightarrow H^*(H) {\buildrel {\times p} \over \longrightarrow} H^*(H) {\buildrel red_{p^*} 
      \over \longrightarrow} H^*(H;{\mathbb F}_p) {\buildrel \delta \over \longrightarrow} \ldots$}
  \put(4,1.7){\vector(0,-1){1.2}}
  \put(6,1.7){\vector(0,-1){1.2}}
  \put(8,1.7){\vector(0,-1){1.2}}
  \put(4.2,1){$res^G_H$}
  \put(6.2,1){$res^G_H$}
  \put(8.2,1){$res^G_H$}
 \end{picture}

\section{Initial Calculations}
We determine the integral cohomology of the five finite groups that form the basis of the Bianchi group's cohomology, then apply techniques from the last section to determine some intermediate results.  Subscripts on classes denote the degree of the class.  

\begin{lemma} $H^*(\Z/2;\Z) \cong \Z[x_2] / (2x_2 = 0)$ and $H^*(\Z/3;\Z) \cong \Z[x_2] / (3x_2 = 0)$.
\end{lemma}
  
\begin{lemma}
The 2-primary component $H^*(\St;\Z)_{(2)} \cong H^*(\Z/2;\Z) \cong \Z[x_2] / (2x_2 = 0)$, and the isomorphism is induced by restriction.  The 3-primary component $H^*(\St;\Z)_{(3)} \cong \Z[x_4] / (3x_4 = 0)$, and this is the image of $res^{\St}_{\Z/3}$. 
\label{St}
\end{lemma}

\prf Consider the short exact sequence 
$$0 \longrightarrow \Z/3 \longrightarrow \St \longrightarrow \Z/2  \longrightarrow 0.$$
There is a spectral sequence with $E_2^{p,q} \cong H^p(\Z/2; H^q(\Z/3;\Z))$ converging to $H^*(\St;\Z)$.  The horizontal edge of the spectral sequence is isomorphic to $H^*(\Z/2;\Z)$.   The vertical edge, however, is the invariants $H^*(\Z/3;\Z)^{\Z/2}$.  The horizontal edge yields the 2-primary part of the cohomology.  On the vertical edge, the action is free in dimensions $2 + 4n$ and trivial otherwise \cite{Bro}.  \fin

We use the Bockstein spectral sequence to determine the integral cohomology of $\Dt$ and $\Af$.

\begin{lemma}
$H^*(\Dt; \Z) \cong \Z[y_2, z_2, y_3]$ with relation $y_3^2 = y_2^2 z_2 + y_2 z_2^2$.  All classes have order 2.
\end{lemma}

\prf Let $x_1$ and $y_1$ generate $H^1(\Dt; \F)$.  Since $Sq^1(x_1) = x_1^2$ and $Sq^1(y_1) = y_1^2$, there are classes $y_2$ and $z_2 \in H^2(\Dt; \Z)$ with $red_{2^*}(y_2) = x_1^2$ and $red_{2^*}(z_2) = y_1^2$.  We will identify $y_2$ and $z_2$ with their images under $red_{2^*}$.  This directly implies that $y_2$ and $z_2$ generate a rank-2 polynomial algebra over $\Z$, as there are no relations between $x_1^2$ and $y_1^2$.  There is another non-trivial square, $Sq^1(x_1 y_1) = x_1^2 y_1 + x_1 y_1^2$, which represents the generator $y_3 \in H^3(\Dt; \Z)$.  Since $0 \neq (x_1^2 y_1 + x_1 y_1^2)^n \in H^*(\Dt; \F)$, $y_3$ is also a polynomial class.  Furthermore, as 
$$ (x_1^2 y_1 + x_1 y_1^2)^2 = x_1^2 y_1^2 (x_1^2 + y_1^2)$$
we get the relation $y_3^2 = y_2 z_2 (y_2 + z_2)$ in the integral cohomology.  To finish the argument, note that all classes whose square are 0 are also in the image of $Sq^1$.  Thus, the spectral sequence collapses at the $E_1$ page, and all classes in the integral cohomology have order 2.  A Poincar\'e series confirms that we have found all the relations. 
\fin

\begin{lemma}
$H^*(\Af; \Z)_{(2)} \cong \Z[y_3,y_4,y_6]$ with relation $y_3^4 + y_4^3 + y_3^2 y_6 + y_6^2 = 0$.  All classes have order 2.  $H^*(\Af; \Z)_{(3)} \cong H^*(\Z/3 ;\Z)$, and this isomorphism is given by $res^{\Af}_{\Z/3}$.  \label{Af}
\end{lemma}

\prf  At the prime 2, start with $H^*(\Af;\F) \cong \F[u_2,v_3,w_3]$ with relation $u_2^3 + v_3^2 + w_3^2 + v_3w_3=0$ (see \cite{A-M} for a derivation).  The Bocksteins are $Sq^1(u_2) = v_3$, $Sq^1(v_3) = 0$, and  $Sq^1(w_3) = u_2^2$.  We identify three classes in $H^*(\Af;\Z)$ in low degrees: $y_3$, whose image is $v_3 = Sq^1(u_2)$; $y_4$, whose image is $u_2^2=Sq^1(w_3)$; and $y_6$, whose image is $u_2^3 + v_3w_3 = Sq^1(u_2w_3)$.  $Sq^1(u_2v_3) = v_3^2$ is also a non-trivial square, but this class is the image of $y_3^2$ and has thus already been accounted for.  

One easily confirms that the classes $v_3$, $u_2^2$, $v_3^2$ and $u_2^3 + v_3w_3$ generate the kernel of $Sq^1$ up to degree 6; therefore, $y_3$, $y_4$, and $y_6$ are of order 2.  A Poincar\'e series argument derived from the Universal Coefficient Theorem implies that these classes satisfy a single relation in degree 12.  A straightforward search yields the one above. 

At the prime 3, we note that $\Z/3$ is self-normalizing in $\Af$.  $H^*(\Af; \Z)_{(3)}$ maps isomorphically onto the invariants $H^*(\Z/3; \Z)^{\Af}$ (\cite{Bro}, III.10.3), which by direct calculation is $H^*(\Z/3; \Z)$.  The invariants are also the image of the restriction map.  
\fin

Most of the restriction maps we will need are straightforward.  However, $res^{\Af}_{\Z/2}$ is somewhat more complicated.

\begin{lemma}
The map $res^{\Af}_{\Z/2}$ on integral cohomology sends $y_4^i y_6^j$ to $x_2^{2i+3j}$ and all other classes to $0$,
where $x_2$ generates $H^*(\Z/2;\Z)$.  \label{AftoZ2}
\end{lemma}

\prf We use the commutative square:

\bigskip
\setlength{\unitlength}{1cm}
\begin{picture}(8,2.6)
  \put(3.5,2){$ H^*(\Af ;\F) {\buildrel \delta \over \longrightarrow} H^{*+1}(\Af; \Z)$}
  \put(3.5,0){$ H^*(\Z/2 ;\F) {\buildrel \delta \over \longrightarrow} H^{*+1}(\Z/2; \Z)$}
  \put(4.6,1.7){\vector(0,-1){1.2}}
  \put(7.6,1.7){\vector(0,-1){1.2}}
  \put(4.8,1){$res^{\Af}_{\Z/2}$}
  \put(7.8,1){$res^{\Af}_{\Z/2}$}
\end{picture}

\medskip

Let $x_1$ and $x_2$ generate $H^1(\Z/2;\F)$ and $H^2(\Z/2;\Z)$ respectively.  From results in \cite{A-M}, $res^{\Af}_{\Z/2}$ on the left sends $u_2$ to $x_1^2$, $v_3$ to $0$, and $w_3$ to $x_1^3$.  By extension, $u_2^i w_3^j$ is mapped to $x_1^{2i+3j}$ if $j$ is odd; all other classes are mapped to zero.  The bottom map sends $x_1^{2i-1}$ to $x_2^i$.  Since we determined the top map in Lemma \ref{Af}, we should be able to determine the restriction map on the right.  

Consider the class $u_2^{2i+1} w_3^{2j-1}$, which maps non-trivially to $x_2^{2i+3j}$ under the composition of the left and bottom maps.   $Sq^1(u_2^{2i+1} w_3^{2j-1}) = u_2^{2i} w_3^{2j-2} (u_2^3 + v_3 w_3)$, which corresponds to the class $y_4^i (y_3^2 + y_6)^{j-1} y_6$ in the integral cohomology.  We claim that the only component in the expansion of this product with a non-trivial image under $res^{\Af}_{\Z/2}$ is the term $y_4^i y_6^j$.

Note that if $Sq^1(z) = 0$, then $Sq^1 (w z) = (Sq^1 w)z$ for any class $w$.  In $H^*(\Af;\F)$, $u_2^3 + v_3 w_3 = v_3^2+w_3^2$.  Thus, $Sq^1 (u_2^3 + v_3 w_3) = 0$, so  $Sq^1 \left( u_2^{2j+1} v_3^{i-1} (u_2^3 + v_3 w_3)^k \right)  = u_2^{2j} v_3^i (u_2^3 + v_3 w_3)^k$, which corresponds to the class $y_3^i y_4^j y_6^k$.  On the other hand, since any term in $H^*(\Af;\F)$ containing a $v_3$ maps to zero under $res^{\Af}_{\Z/2}$, $y_3^i y_4^j y_6^k$ must map to $0$ also.    

To finish the argument, $Sq^1 (u_2^{2i-2} w_3)$ corresponds to $y_4^i$ and maps non-trivially onto $x_2^{2i}$.  Therefore, the restriction map is only non-zero when it sends $y_4^i y_6^j$ to $x_2^{2i+3j}$. \fin

\begin{remark}  Using a similar, but easier, analysis, one can show that $res^{\Dt}_{\Z/2}$ sends powers of $x_2 \in H^2(\Dt; \Z)$ isomorphically to powers of the corresponding generator of $H^2(\Z/2;\Z)$ and vanishes otherwise. \label{DtontoZ2}  
\end{remark}

We next move to the calculation of three amalgamated products that appear as base groups for a number of the Bianchi groups presented in this article.  In many cases, a Bianchi group's cohomology is nearly identical to its base group's, so the following results consolidate later calculations.  

\begin{lemma} $H^*(\Dt \ast_{\Z/2} \Dt ;\Z)_{(2)} \cong \Z[x_2, y_2, z_2, y_3, z_3]$ with relations $y_2 z_3 = z_2 y_3 = 0$, $y_3^2 = x_2 y_2 (x_2 + y_2)$, and ${z_3}^2 = x_2 z_2 (x_2 + z_2)$.  All classes have order 2. \label{DtZtDt}
\end{lemma}

\prf We use Formula (\ref{LES}) with integer coefficients:
$$ \longrightarrow H^*(\Dt \ast_{\Z/2} \Dt;\Z) \brarrow H^*(\Dt;\Z) \oplus H^*(\Dt;\Z) \ararrow H^*(\Z/2;\Z) \longrightarrow $$
By Remark \ref{DtontoZ2} $\alpha$ is surjective, so this sequence breaks into short exact sequences.  Let $x_2,\ y_2,\ y_3$ and $x'_2,\ z_2,\ z_3$ be generators in each copy of $H^*(\Dt;\Z)$ in the direct sum, and let $w_2$ be the generator of $H^2(\Z/2)$ that is the restriction of both $x_2$ and $x'_2$ under $\alpha$.  Now $y_2,\ z_2,\ y_3,\ z_3$ and $x_2 + x'_2$ are all in $ker (\alpha)$ and have order 2.  Since this kernel is isomorphic to the image of $\rho$, we can construct the cohomology ring structure of $\Dt \ast_{\Z/2} \Dt$ at this point.  We abuse notation and refer to classes in $H^*(\Dt \ast_{\Z/2} \Dt;\Z)$ by their image under $\rho$, although we rename $x_2 + x'_2$ as $\bar x_2$.  Most of the relations are obvious; note that $y_3^2 = \bar x_2 y_2 (\bar x_2 + y_2)$ as products of $x'_2$ with $y_2$ equal $0$.  That the relations we found form a minimal set follows from a Poincar\'e series argument.  We drop the bar on the $x_2$ for brevity.

One can derive the same result by using the Bockstein spectral sequence.  Briefly, start with $H^*(\Dt \ast_{\Z/2} \Dt ;\F) \cong \F [p_1, q_1, r_1] /p_1 q_1 = 0$ \cite{Be} and make the following identifications: $x_2 \leftrightarrow r_1^2, y_2 \leftrightarrow p_1^2, z_2 \leftrightarrow q_1^2, y_3 \leftrightarrow p_1 r_1(p_1 + r_1)$, and  $z_3 \leftrightarrow q_1 r_1 (q_1 + r_1)$.  These are the Bocksteins of the classes $r_1, p_1, q_1, p_1 r_1$, and $q_1 r_1$ respectively.  It is straightforward to confirm that the stated relations hold.  \fin

For $\Dt \ast_{\Z/2} \Dt$, the integral calculations are possible using either the Bockstein spectral sequence or Formula \ref{LES}.  In practice, both techniques are often needed to determine the answer.  The long exact sequence is an effective means to determine generators in low degree, and show exactly how classes arise from pieces of the amalgamated product.  The Bockstein spectral sequence, on the other hand, is helpful in the determination of the relations between generators.  

\begin{lemma}
$H^*(\Af \ast_{\Z/2} \Dt ; \Z)_{(2)} \cong \Z[y_2, y_3, z_3, y_4, z_4, y_5, y_6]$.  All classes are of order two, and satisfy the relations given in the table below. \label{AfZtDt}
\end{lemma}

\prf By results in \cite{Be}, $H^*(\Af \ast_{\Z/2} \Dt;\F) \cong \F [x_1,y_2,u_2,v_3,w_3]$ with relations $u_2^3 + v_3^2 + w_3^2 + v_3w_3 = 0,\ x_1^2u_2 = y_2^2,\ x_1w_3 = u_2y_2,\ y_2w_3 = x_1u_2^2$, and $x_1v_3 =
y_2v_3 = 0$.  In addition, $Sq^1 (x_1) = x_1^2$, $Sq^1 (y_2) = x_1 (u_2 + y_2)$, $Sq^1 (u_2) = v_3$, and $Sq^1 (w_3) = u_2^2$.  

We apply the Bockstein spectral sequence and identify generators for the integral cohomology ring.  In general, one writes down classes in the mod-2 cohomology by degree and determines their Steenrod squares.  Every time an integral class is encountered that is not in the ideal generated by the previous classes, it is added to the generating set.  Relations among classes of order 2 can be determined using the relations in the mod-2 cohomology ring.  There are some ways to speed up this occasionally difficult process.  Once all integral generators have been found, it is straightforward to write down all possible combinations of generators by degree.  Using the Poincar\'e series for the mod-2 cohomology and the Universal Coefficient Theorem, one knows how many classes are in the integral cohomology in each degree; this determines the number of relations among the classes.  One continues this process until an integral Poincar\'e series argument or other appropriate technique shows that all generators and relations have been found.    

Following this process, we identify generators $y_2 \leftrightarrow x_1^2 = Sq^1(x_1);\ y_3 \leftrightarrow x_1 (u_2 + y_2) = Sq^1(y_2);\ z_3 \leftrightarrow v_3 = Sq^1(u_2);\ y_4 \leftrightarrow x_1^2 u_2 = Sq^1(x_1 u_2);\ z_4 \leftrightarrow u_2^2 = Sq^1(w_3);\ y_5 \leftrightarrow x_1^2 w_3 + x_1 u_2^2 = Sq^1(x_1 w_3)$;\ and $y_6 \leftrightarrow u_2^3 + v_3 w_3 = Sq^1(u_2 w_3).$   In this example, the number of generators makes it difficult to confirm that the spectral sequence collapses at $E_2$.  However, using Formula \ref{LES} we get the sequence
$$ \rightarrow H^*(\Af \ast_{\Z/2} \Dt;\Z)_{(2)} \brarrow H^*(\Af;\Z)_{(2)} \oplus H^*(\Dt;\Z) \ararrow H^*(\Z/2;\Z) \rightarrow .$$
The map $\alpha$ is a surjection, which follows from Lemma \ref{AftoZ2} and Remark \ref{DtontoZ2}.  Therefore, the long exact sequence breaks into short exact sequences and all 2-torsion in $H^*(\Af \ast_{\Z/2} \Dt;\Z)_{(2)}$ must be of order 2.  This implies that $E_1 = E_{\infty}$.  We include a table of relations in the cohomology ring, organized by degree.

\bigskip

\begin{center}
\begin{tabular}{|l|l|}
  \hline
  Degree & Relation(s) \\
  \hline
  5      & $y_2z_3=0$ \\
  \hline
  6      & $y_3 z_3 = 0$, $y_3^2 = y_2(y_4+z_4)$ \\
  \hline
  7      & $z_3 y_4 = 0$, $y_3 y_4 = y_2 y_5$ \\
  \hline
  8      & $z_3 y_5 = 0$ , $y_4^2 = y_2^2 z_4$, $y_4 z_4 = y_2 y_6$, $y_3 y_5 = y_4^2 + y_2 y_6$\\
  \hline
  9      & $y_3 y_6 = z_4 y_5$, $y_2 y_3 z_4 = y_4 y_5$ \\
  \hline
  10     & $y_5^2 = y_3^2 z_4$, $y_4 y_6 = y_2 z_4^2$ \\
  \hline
  11     & $y_5 y_6 = y_3 z_4^2$ \\
  \hline
  12     & $y_6^2 = y_6 z_3^2 + z_3^4 +z_4^3$ \\
  \hline
\end{tabular}  
\end{center}
As a sample confirmation of the relations, note that $y_5^2 \leftrightarrow x_1^2 u_2^4 + x_1^4 w_3^2$ and $y_3 z_4^2 \leftrightarrow x_1^2 u_2^4 + x_1^2 u_2^2 y_2^2$.  Using the relation $x_1 w_3 = u_2 y_2$ from the mod-2 ring, we find that the two expressions are equal.   \fin

\begin{lemma}
$H^*(\Af \ast_{\Z/2} \Af ; \Z)_{(2)} \cong \Z[y_3, y_4, y_6, z_6, y_7, y_9] \oplus \Z[t_3]$.  All classes have order 2 with the exception of $t_3$, which has order 4.  The relations are given in the table below. \label{AfZtAf}
\end{lemma}

\prf $H^*(\Af \ast_{\Z/2} \Af;\F) \cong \F[u_2,v_3, \bar v_3,w_3] \oplus \F[s_2]$ with relations $u_2^3 + w_3^2 + v_3^2 + \bar v_3^2 + w_3 (v_3 + \bar v_3) = 0,\ v_3 \bar v_3 = 0$ \cite{Be}.  Also, $Sq^1 (u_2) = v_3 + \bar v_3$ and $Sq^1 (w_3) = u_2^2$.

We identify the following classes in the integral cohomology ring: $y_3 \leftrightarrow v_3 + \bar v_3 = Sq^1(u_2);\ y_4 \leftrightarrow u_2^2 = Sq^1(w_3);\ y_6 \leftrightarrow u_2^3 + (v_3 + \bar v_3) w_3 = Sq^1(u_2 w_3);\ z_6 \leftrightarrow v_3^2 = Sq^1(u_2 v_3);\ y_7 \leftrightarrow u_2^2 v_3 = Sq^1(v_3 w_3);$ \ and $y_9 \leftrightarrow v_3 (u_2^3 + v_3 w_3) = Sq^1(u_2 v_3 w_3)$.  An exhaustive search finds all of the relations among these classes. 
\bigskip

\begin{center}
\begin{tabular}{|l|l|}
  \hline
  Degree & Relation(s) \\
  \hline
  10      & $y_3 y_7 = y_4 z_6$ \\
  \hline
  12      & $y_3^4 + y_4^3 + y_3^2 y_6 + y_6^2 =0$, $y_3 y_9 = y_6 z_6$, $z_6^2 = y_3^2 z_6$\\
  \hline
  13      & $y_6 y_7 = y_4 y_9$, $z_6 y_7 = y_3^2 y_7$ \\
  \hline
  14      & $y_7^2 = y_4^2 z_6$ \\
  \hline
  15      & $z_6 y_9 = y_3 y_6 z_6$ \\
  \hline
  16      & $y_7 y_9 = y_4 y_6 z_6$ \\
  \hline
  18      & $y_9^2 = y_3 y_6 y_9$, $y_9^2 + y_3 z_6 y_9 + y_4 y_7^2 + z_6^3 = 0$ \\
  \hline
\end{tabular}  
\end{center}

Note that $Sq^1 (v_3) = Sq^1 \bar v_3 = 0$, whereas $Sq^1 (u_2) = v_3 + \bar v_3$.  That means that there is a class in the integral cohomology which has at least order 4.  We show that this class has exactly order 4 by considering the short exact sequence of 2-primary groups.
$$ 0 \ararrow H^2(\Z/2;\Z) \drarrow H^3(\Af \ast_{\Z/2} \Af;\Z)_{(2)} \brarrow H^3(\Af;\Z)_{(2)} \oplus H^3(\Af;\Z)_{(2)} \ararrow 0.$$
Using Lemma \ref{Af}, we substitute to get 
$$ 0 \longrightarrow \Z/2 \longrightarrow H^3(\Af \ast_{\Z/2} \Af;\Z)_{(2)} \longrightarrow \Z/2 \oplus \Z/2 \longrightarrow 0.$$
We know from the Bockstein spectral sequence that there are two classes in $H^3(\Af \ast_{\Z/2} \Af;\Z)_{(2)}$, exactly one of which has order greater than two.  This implies that $H^3(\Af \ast_{\Z/2} \Af;\Z)_{(2)} \cong \Z/2 \oplus \Z/4$.  We denote the class of order 4 by $t_3$.  

It remains to determine products with $t_3$.  First, $t_3^2 = 0$, as $t_3$ is the only class in the cohomology ring of order 4.  Without loss of generality, we can say that $\rho (t_3) = v_3$.  Since $\rho (x_3) = v_3 + \bar v_3$, if $t_3 x_3$ is a non-zero product it must be $x_6$, as both classes map to $v_3^2$ under $\rho$ and there is only one such class.   But $x_6^2 \neq 0$, whereas $(t_3 x_3)^2 = 0$; we conclude $t_3 x_3 = 0$.  Similar arguments show that all products with $t_3$ are 0.  \fin

\section{The Bianchi Groups $\G2$ and $\G6$}
In this section, we calculate the integral cohomology rings of $\G2$ and $\G6$ in detail.  We use the Bockstein spectral sequence to determine the 2-torsion, and Formula \ref{LES} for the rest of the ring structure.  In here and what follows, Roman letters refer to torsion classes, and Greek letters to torsion-free classes.

\begin{theorem}
$H^*(\Gamma_2;\Z)_{(2)} \cong \Z[y_2, t_3, y_3, y_4, y_5, y_6, y_7, \sigma_1]$.  All torsion classes have order 2, except for $t_3$ which has order 4.   These classes satisfy the relations given in the table below, and $\sigma_1 y_2 = 2 t_3$; all other products with $\sigma$ and $t_3$ are trivial.  $H^*(\G2;\Z)_{(3)} \cong \Z[x_2, \sigma_1]$ with relations $\sigma_1^2 = 3x_2 = 0$.
\end{theorem}

\prf  From \cite{Be}, $H^*(\Gamma_2;\F) \cong \F[n_1,m_2,n_3, m_3] \oplus \F[s_1] \oplus  \F[s_2]$ with relations $n_1n_3 = 0$ and  $m_2^3 + m_3^2 + n_3^2 + m_3n_3 + n_1m_2m_3 = 0$. All products with ${\sigma}_1$ vanish except for the product ${\sigma}_1 n_1 = \sigma_2$.  The non-trivial Bocksteins are $Sq^1 (n_1) = n_1^2$, $Sq^1 (m_2) = n_1 m_2 + n_3$, and $Sq^1 (m_3) = m_2^2$.

We apply the Bockstein spectral sequence and identify generators $y_2 \leftrightarrow n_1^2 = Sq^1(n_1);\ y_3 \leftrightarrow n_1 m_2 + n_3 = Sq^1(m_2);\ y_4 \leftrightarrow m_2^2 = Sq^1(m_3);\ y_5 \leftrightarrow n_1^2 m_3 + n_1 m_2^2 = Sq^1(n_1 m_3);\ y_6 \leftrightarrow n_1 m_2 m_3 + n_3 m_3 + m_2^3 = Sq^1(m_2 m_3)$; and $y_7 \leftrightarrow m_2^2 n_2 = Sq^1(m_3 n_3)$.   The following table contains the relations by degree. 

\bigskip

\begin{center}
\begin{tabular}{|l|l|}
  \hline
  Degree & Relation(s) \\
  \hline
  8       & $y_2^2 y_4 = y_2 y_3^2$, $y_2 y_6 = y_3 y_5$ \\
  \hline
  9       & $y_2 y_7 = 0$ \\
  \hline
  10      & $y_2 y_4^2 = y_3 y_7 + y_3^2 y_4$, $y_5^2 = y_2^2 y_6 + y_2 y_4^2$  \\
  \hline
  11      & $y_5 y_6 = y_3^2 y_5 + y_3 y_4^2 + y_4 y_7$, $y_2 y_3 y_6 = y_2 y_4 y_5$ \\
  \hline
  12      & $y_5 y_7 = 0$, $y_6^2 = y_3^2 y_6 + y_3^4 + y_2^2 y_4^2 + y_4^3$ \\
  \hline
  13      & $y_6 y_7 = y_3 y_4 y_6 + y_4^2 y_5$ \\
  \hline
  14      & $y_7^2 = y_2 y_4^3 + y_3^2 y_4^2$ \\
  \hline
\end{tabular}  
\end{center}
It remains to determine what happens to the classes $s_1$ and $s_2$, and whether there are classes in the integral cohomology with order 4 or higher.  We use Fl\"oge's presentation for the group:
$\Gamma_2 =\ <A,V,S,M,U; A^2 = S^3 = (AM)^2 = M^2 = V^3 = 1,\ AM = SV^2,\ U^{-1}AU = M,\ 
U^{-1}SU = V>$ \cite{Fl}.  Set $G =\ <A,V,S,M; A^2 = M^2 = (AM)^2 = S^3 = V^3 = 1,\ AM = SV^2>$ and consider the subgroups
\begin{eqnarray*}
  G_1 & = & <A,M; A^2 = M^2 = (AM)^2 = 1>\   \cong \Dt , \\
  G_2 & = & <S,V; S^3 = V^3 = (SV^2)^2 = 1>\ \cong \Af , \\
  H   & = & <AM = SV^2,\ (AM)^2 = 1 >\       \cong \Z/2 . 
\end{eqnarray*}
Let $G$ be the amalgamated product $G_1 {\ast}_H G_2 \cong \Dt {\ast}_{\Z/2} \Af$.  Then $\Gamma_2 = G \ast_{\Ps}$, where $\Ps =\ <A,S>$ and the twisting is induced by the group element $U$. 

We have already calculated $H^*(G;\Z)_{(2)}$ in Lemma \ref{AfZtDt}, but in order to keep track of the classes more carefully, we look at this example again from the point of view of group structure.  The amalgamated subgroup in $G$ is $<AM>$, so let the two classes $u_2, v_2 \in H^2(\Dt;\Z)$ correspond to the subgroups $<AM>$ and $<A>$.  Note that $\alpha$ sends $u_2$ to the generator of $H^2(\Z/2;\Z)$ but sends $v_2$ to 0.  The map $\alpha$ is a surjection by Lemma \ref{AftoZ2} and Remark \ref{DtontoZ2}, so the long exact sequence splits into short exact sequences.  It follows that in $H^*(G;\Z)_{(2)}$ there is a single class in degree 2 and two classes in degree 3, all of order 2.  

Next consider the HNN extension that generates $\G2$, where the group element $U$ twists $<A>$ to $<M>$.  As $\Ps \cong \Z/2 \ast \Z/3$, let $w_2$ be the generator of $H^2(\Ps;\Z)$ corresponding to the subgroup $<A>$.  Then $res^{\G2}_{\Ps} (v_2) = \theta^*(v_2)= w_2$, so $\alpha(v_2) = 0$.  We can fit this information into Formula \ref{LES} for an HNN extension:
$$ 0 \rightarrow H^2(\Gamma_2) \brarrow \Z/2 \ararrow \Z/2 \drarrow H^3(\Gamma_2) \brarrow \Z/2 \oplus \Z/2 \rightarrow 0.$$
The map $\alpha$ is 0, so $\rho$ is an isomorphism.  We abuse notation and denote the class in $H^2(\G2;\Z)_{(2)}$ by $v_2$ as well.  (Notice that $res^{\G2}_{\Z/2} (v_2) = w_2$, as both classes are dual to the group element $A$.)  Since the Bockstein spectral sequence guarantees the existence of a single class of order 2 in $H^3(\Gamma_2)$, there must be another class in degree 3 of order 4, denoted by $t_3$.  Also, there is a torsion-free HNN class in degree 1, $\sigma_1$.  Now by Lemma \ref{compat}, $\sigma_1 v_2 = \delta(1) \cup v_2 = \delta(1 \cup res^{\G2}_{\Z/2}(v_2)) = \delta(w_2) = 2 t_3$.  Finally, one can show that $\alpha$ is a surjection in degrees 3 and higher.  This implies that $t_3$ is the only 2-torsion class of order greater than 2, which finishes the argument.

To determine $H^*(\G2;\Z)_{(3)}$, we have that $H^*(G;\Z)_{(3)} \cong H^*(\Z/3;\Z)$, so all that remains is to determine the effect of the final HNN extension.  The group element $U$ sends $S$ to $V$, which are both elements of degree 3 in $\Af$.  Therefore, both the injection and twisting maps are isomorphisms, which implies that $\alpha$ is the zero map.  We get short exact sequences
$$ 0 \rightarrow H^*(\Ps;\Z)_{(3)} \drarrow H^{*+1}(\G2;\Z)_{(3)} \brarrow H^{*+1}(\Z/3;\Z)_{(3)} \ararrow 0$$
This yields classes $x_2^n \in H^{2n}(\G2;\Z)_{(3)}$ and classes $\delta(w_2^n) \in H^{2n+1}(\G2;\Z)_{(3)}$, where $w_2$ is the generator of $H^2(\Ps;\Z)_{(3)}$.  All classes have order 3.  There is also a torsion-free HNN class, $\delta(1) \in H^1(\G2;\Z)_{(3)}$, which is the class $\sigma_1$ we identified before.  We claim that products with this class generate the odd dimensional cohomology.  For note that $res^{\G2}_{\Ps} (x_2^n) = w_2^n$.  Thus, $\delta(w_2^n) = \delta(w_2^n \cup 1) = \delta(res^{\G2}_{\Ps} (x_2^n) \cup 1) = x_2^n \cup \delta(1) = x_2^n \cup \sigma_1$.  
\fin

\begin{theorem}  $H^*(\Gamma_6 ; \Z)_{(2)} \cong \Z[y_3, y_4, y_6, z_6, y_7, y_9] \oplus \Z[y_2, \tau_1] \oplus \Z[t_3] \oplus \Z[\sigma_1] \oplus \Z[\sigma_2]$.  All torsion classes have order 2 except for $t_3$, which is of order 4.  Relations among these classes are almost identical to $H^*(\Af \ast_{\Z/2} \Af; \Z)_{(2)}$.   Products with torsion-free classes are trivial except for $y_2^k \tau_1$.  At the prime 3, $H^*(\G6;\Z)_{(3)} \cong \Z[x_2, \tau_1] \oplus \Z[\sigma_1] \oplus \Z[\sigma_2]$ with relations $\tau_1^2 = \sigma_1^2 = \sigma_2^2 = 3 x_2 = 0$. 
\end{theorem}

\prf $H^*(\Gamma_6;\F) \cong \F[y_1, u_2,v_3, \bar v_3,w_3] \oplus \F[s_1] \oplus \F[t_1] \oplus \F[s_2] \oplus \F[t_2]$ with relations $u_2^3 + w_3^2 + v_3^2 + \bar v_3^2 + w_3 (v_3 + \bar v_3) = 0,\ v_3 \bar v_3 = 0$, and all products of other classes with each other 0 except for $y_1^j t_1$.   The Bocksteins are $Sq^1(y_1) = y_1^2$, $Sq^1(u_2) = v_3 + \bar v_3$ and $Sq^1(w_3) = u_2^2$ \cite{Be}.

As $\G6$ and $\Af \ast_{\Z/2} \Af$ have almost identical mod-2 cohomology, it is no surprise that most of their generators and relations in the integral cohomology are identical too;  we refer the reader to Lemma \ref{AfZtAf} for the definitions of $y_3, y_4, y_6, z_6, z_7$, and $y_9$, and the table of their relations.  However, $\G6$ has a new class $y_2$ which corresponds to $y_1^2 = Sq^1(y_1)$.  We switch techniques and directly calculate the cohomology of $\G6$ using its group structure to find relations with $y_2$.

Fl\"oge's presentation for $\Gamma_6$ is $<A,B,M,R,S,U,W; A^2=B^2=M^2=R^3=S^3=(BR)^3=(BS)^3=1,\ AS=MR,\ U^{-1}AU = M,\ U^{-1}SU=R,\ W^{-1}MW=A,\ W^{-1}RBW=SB>$ \cite{Fl}. Let $G$ be the subgroup generated by the group elements $A,B,M,R,S$ and their relations.  The following subgroups of $G$, 
 \begin{eqnarray*}
  G_1 & = & <S,A,B; A^2=S^3=B^2=(BS)^3=1 > \ \cong \Z/2 \ast \Af,   \\
  G_2 & = & <R,B,M; M^2=R^3=B^2=(BR)^3=1 > \ \cong \Z/2 \ast \Af,   \\
  H   & = & <B,AS=MR; B^2 =1 >             \ \cong \Z/2 \ast {\Z}, 
\end{eqnarray*}
combine to form $G \cong G_1 {\ast}_H G_2 \cong (\Z/2 \ast \Af) {\ast}_{(\Z/2 \ast {\mathbb Z}) } (\Z/2 \ast \Af)$.  The associated subgroups of the HNN extensions, $<A,\ S>$ and $<M,\ RB>$, are both isomorphic to $\Ps$.  We add the group element $U$ to get the HNN extension $G_3 = G \ast_{\Ps}$.  Finally, we add $W$ to get $\Gamma_6 = G_3 \ast_{\Ps}$.  Thus, the cohomology of $\Gamma_6$ can be calculated in three steps.

We use Formula \ref{LES} to find $H^*(G;\Z)_{(2)}$.  
$$ \ldots \drarrow H^*(G;\Z)_{(2)} \brarrow H^*(\Z/2 \ast \Af;\Z)_{(2)} \oplus H^*(\Z/2 \ast \Af;\Z)_{(2)} \ararrow H^*(\Z/2 \ast \Z;\Z)  \drarrow \ldots.$$

From the group presentation, we see that the amalgamated copy of $\Z/2$ injects into both copies of $\Af$, forming $\Af \ast_{\Z/2} \Af$ as a subgroup.  The only torsion-free class in the sequence, in $H^1(\Z;\Z)$, is not hit by $\alpha$.  Its image under $\delta$ yields a torsion-free class in degree 2 which we denote by $\sigma_2$.  Since this class and $2 t_3$ are the only ones in the image of $\delta$, they have trivial products with the rest of the ring elements by Lemma \ref{compat} and with each other by Lemma \ref{PY}.  The other copies of $\Z/2$, corresponding to the group elements $A$ and $M$, go to 0 under $\alpha$.  Thus, they appear as two copies of $H^*(\Z/2)$ in $H^*(G;\Z)_{(2)}$, which is isomorphic to the reduced direct sum $H^*(\Af \ast_{\Z/2} \Af;\Z)_{(2)} \widetilde{\oplus} H^*(\Z/2;\Z) \widetilde{\oplus} H^*(\Z/2;\Z) \oplus \Z[\sigma_2]$.  Let $v_2$ and $w_2$, dual to $<A>$ and $<M>$, generate the two copies of $H^2(\Z/2;\Z)$. 

For the first HNN extension, the element $U$ twists $A$ to $M$, with $<A> \subseteq \Ps$, the associated subgroup of the extension.  By Formula \ref{LES}, we get the sequence 
$$ \ldots \drarrow H^*(G_3;\Z)_{(2)} \brarrow H^*(G;\Z)_{(2)} \ararrow H^*(\Ps;\Z)_{(2)}  \drarrow \ldots.$$
Now as $\alpha$ is a surjection in positive degrees the long exact sequence breaks into short exact sequences.
In $H^*(G;\Z)_{(2)}$, the sum of the two generators $v_2 + w_2$ is in the kernel of $\alpha$ because of the twisting; let $y_2$ be the inverse image under $\rho$ of this new generator.  There is also an HNN class, $\sigma_1$, which has trivial products as it is the only class that arises from the image of $\delta$.  Thus, $H^*(G_3;\Z)_{(2)} \cong H^*(\Af \ast_{\Z/2} \Af;\Z)_{(2)} \oplus H^*(\Z/2) \oplus \Z[\sigma_1] \oplus \Z[\sigma_2]$.  

The final HNN extension yields $H^*(\Gamma_6;\Z)_{(2)}.$  Here the twisting sends $M$ back to $A$.  In terms of the maps induced by inclusion of twisting, $res^{\G6}_{\Ps} (w_2) = \theta^*(v_2)$, so $\alpha(y_2) = 0$.  We abuse notation and let $y_2$ also denote its inverse image in $H^*(\G6;\Z)_{(2)}$ that arises from the kernel of $\alpha$.  The long exact sequence of Formula \ref{LES} breaks into short exact sequences that include a degree shift:
$$ 0 \ararrow H^*(\Ps;\Z)_{(2)} \drarrow H^{*+1}(\Gamma_6;\Z)_{(2)} \brarrow H^{*+1}(G_3;\Z)_{(2)} \ararrow 0.$$
New classes in $H^*(\Gamma;\Z)_{(2)}$ include the HNN class $\tau_1 = \delta(1)$ and $\delta(u_2^j)$ in all odd degrees, where $u_2 \in H^2(\Ps;\Z)$ corresponds to the copy of $\Z/2$ generated by the group element $M$.  We use Lemma \ref{compat} to show that the new classes in odd degree are actually products of other classes.  Let $i$ denote the subgroup inclusion $\Ps \hookrightarrow \G6$.  By construction, $res^{\G6}_{\Ps} (y_2) = u_2$.  Therefore, $\delta(u_2^k) = \delta(res^{\G6}_{\Ps} (y_2^k) \cup 1) = y_2^k \cup \delta(1) = y_2^k \tau_1$.  

The 3-primary cohomology calculation is similar.  First, $H^*(G;\Z)_{(3)} \cong H^*(\Z/3;\Z) \widetilde{\oplus} H^*(\Z/3;\Z) \oplus \Z[\sigma_2]$.  The torsion-free class, $\sigma_2$, is the image of $H^1(\Z;\Z)$ under $\delta$.  The first HNN extension has the effect of identifying the two polynomial classes as it twists one copy of $\Z/3$ to the other.  Thus, $H^*(G_3;\Z)_{(3)} \cong H^*(\Z/3;\Z) \oplus \Z[\sigma_1] \oplus \Z[\sigma_2]$, where $\sigma_1$ is the torsion-free HNN class.

For the final HNN extension, $W$ twists the second copy of $\Z/3$ back to the original one; thus $\alpha$ is 0.  The long exact sequence splits as in the case of $\G2$, with the same results: we get a new exterior HNN class $\tau_1$, which multiplies non-trivially with a new polynomial class in degree 2.  \fin

\section{The Bianchi Groups $\G1$, $\G3$, $\G5$, $\G7$, $\G{10}$, and $\G{11}$  }

\begin{theorem} $H^*(\G1;\Z)_{(2)} \cong H^*(\Af;\Z)_{(2)} \widetilde{\oplus} H^*(\Dt;\Z)_{(2)}$.  $H^*(\G1;\Z)_{(3)} \cong \Z[x_3, x_4]$ with relations $3 x_3 = 3 x_4 = x_3^2 = 0$.
\end{theorem}

\prf Fine shows that $\G1$ has the group presentation $<A, B, C, D; A^3 = B^2 = C^3 = D^2 = (AC)^2 = (AD)^2 = (BD)^2 = (BC)^2 = 1 >$ \cite{F}.  We break this into subgroups.
\begin{eqnarray*}
   G_{11} & = & <A,C; A^3 = C^3 = (AC)^2 = 1>\ \cong \Af  \\
   G_{12} & = & <A,D; A^3 = D^2 = (AD)^2 = 1>\ \cong \St  \\
   G_{21} & = & <B,C; B^2 = C^3 = (BC)^2 = 1>\ \cong \St  \\
   G_{22} & = & <B,D; B^2 = D^2 = (BD)^2 = 1>\ \cong \Dt  
\end{eqnarray*}
Set $G_1 = G_{11} \ast_{<A>} G_{12}$ and $G_2 = G_{21} \ast_{<B>} G_{22}$.  Then $\Gamma_1$ is an amalgamated product $G_1 \ast_{\Ps} G_2$, where $\Ps \cong <C,D>$.

Using Formula \ref{LES} for an amalgamated product, one easily shows that $H^*(G_1;\Z)_{(2)} \cong H^*(\Af;\Z)_{(2)}$, as $\alpha$ maps $H^*(\St;\Z)_{(2)}$ isomorphically onto $H^*(\Z/2;\Z)_{(2)}$.  On the other hand, $H^*(G_2;\Z)_{(2)} \cong H^*(\Dt;\Z)_{(2)} \oplus H^*(\St;\Z)_{(2)}$.  The 2-primary result follows by noting that the copy of $H^*(\St;\Z)_{(2)}$ maps isomorphically onto $H^*(PSL_2(\Z);\Z)_{(2)}$ in the final amalgamated product.  

For the 3-primary component, the restriction maps for $\Af$ and $\St$ to $\Z/3$ (Lemmas \ref{St} and \ref{Af}) imply that both $H^*(G_1 ;\Z)_{(3)}$ and  $H^*(G_2 ;\Z)_{(3)} \cong \Z[x_4]/ (3x_4 = 0)$.  From Formula \ref{LES} there is a long exact sequence 
$$ \ldots H^*(\Gamma_1;\Z)_{(3)} \brarrow H^*(G_1;\Z)_{(3)} \oplus H^*(G_2;\Z)_{(3)} \ararrow H^*(\Ps;\Z)_{(3)} \drarrow \ldots$$
Both classes in degree 4 map onto the square of the generator $w_2 \in H^2(\Ps;\Z)_{(3)}$.  Therefore, we get classes $x_4^n \in H^{4n}(\Gamma_1;\Z)_{(3)}$ from the kernel of $\alpha$ and classes $\delta(w_2^{4n+2}) \in H^{4n+3}(\Gamma_1;\Z)_{(3)}$.  These latter classes are exterior by degree considerations.  Note that $res^{\G1}_{\Ps}(x_4^n) = w_2^{2n}$.  Thus, $\delta(w_2^{2n+1}) = \delta(res^{\G1}_{\Ps}(x_4^n) \cup w_2) = x_4^n \cup \delta(w_2) = x_4^n x_3$ by Lemma \ref{compat}.
\fin

\begin{theorem} $H^*(\G3;\Z)_{(2)} \cong H^*(\Af \ast_{\Z/2} \Af;\Z)_{(2)}$.  $H^*(\G3;\Z)_{(3)} \cong \Z[x_2,x_4]$ with relations $3x_2 = 3x_4 = x_2x_4 = 0$.
\end{theorem}

\prf
From calculations in \cite{Be}, $H^*(\G3; \F) \cong \F[u_2,\ v_3,\ \bar v_3,\ w_3] \oplus \F[s_2]$ with relations $u_2^3 + w_3^2 + v_3^2 + \bar v_3^2 + w_3 (v_3 + \bar v_3) = 0,\ v_3 \bar v_3 = 0$.  This is identical to $H^*(\Af \ast_{\Z/2} \Af; \F)$ given in Lemma \ref{AfZtAf}.

The 3-primary calculations are more involved.  We use the fact that $\G3$ acts on a retract known as the ``Mendoza Complex'' \cite{M} with fundamental domain and isotropy groups \cite{S-V} given in Figure 1:

\begin{figure}[h]
\setlength{\unitlength}{1cm}
\begin{picture}(12,2.5)
  \put(5,0){\line(3,4){1.5}}
  \put(6.5,2){\line(3,-4){1.5}}
  \put(5,0){\line(1,0){3}}
  \put(4.3,0){$\Af$}
  \put(8.3,0){$\Af$}
  \put(6.3,2.2){$\St$}
  \put(6.3,0.3){$\Z/3$}
  \put(5.1,1.1){$\Z/2$}
  \put(7.4,1.1){$\Z/2$}
\end{picture} 
\caption{Fundamental Domain for $\G{3}$}
\label{fig2}
\end{figure}
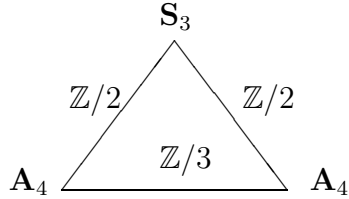

One now uses the equivariant spectral sequence given in \cite{Bro} to show that $H^*(\Gamma_3;\Z)_{(3)} \cong H^*(\Af;\Z)_{(3)} \widetilde{\oplus} H^*(\St;\Z)_{(3)}$.  The details are identical, but easier than the those presented in \cite{Be}.  We refer the interested reader there.
\fin

\begin{theorem} $H^*(\G5 ; \Z)_{(2)} \cong \Z[y_2, z_2, y_3, y_5, z_5, y_6] \oplus \Z[t_3] \oplus \Z[\sigma_1] \oplus \Z[\tau_1] \oplus \Z[\sigma_2]$.  The class $t_3$ has order 4; all other torsion classes have order two.  Relations are given in the table below.  $H^*(\G5;\Z)_{(3)} \cong \Z[x_2, \tau_1] \oplus \Z[\sigma_1] \oplus \Z[\sigma_2]$ with relations $\tau_1^2 = \sigma_1^2 = \sigma_2^2 = 3 x_2 = 0$.

\end{theorem}

\prf From \cite{Be}, the polynomial classes of $H^*(\G5 ; \F)$ generate the ring $\F [l_1, m_1, m_3] / m_3^2 + l_1^2 m_1 m_3 + l_1 m_1^2 m_3 = 0$.  There are also four exterior classes, two each in dimensions 1 and 2, which have trivial products with other classes.  The Bocksteins are $Sq^1 (l_1) = l_1^2$ and $Sq^1 (m_1) = m_1^2$.  Enumerating classes by degree, we find the six generators mentioned above: $y_2 \leftrightarrow l_1^2 = Sq^1 (l_1)$;  $z_2 \leftrightarrow m_1^2 = Sq^1 (m_1)$;   $y_3 \leftrightarrow l_1^2 m_1 + l_1 m_1^2 = Sq^1 (l_1 m_1)$; $y_5 \leftrightarrow l_1^2 m_3 = Sq^1 (l_1 m_3)$; $z_5 \leftrightarrow m_1^2 m_3 = Sq^1 (m_1^2 m_3)$; and $y_6 \leftrightarrow m_3^2 = Sq^1 (l_1 m_1 m_3)$.  This last relation follows as $m_3^2 = l_1^2 m_1 m_3 + l_1 m_1^2 m_3$ in mod-2 cohomology.  These six integral classes satisfy the following relations:

\begin{center}
\begin{tabular}{|l|l|}
  \hline
  Degree & Relation(s) \\
  \hline
  6       & $y_3^2 = y_2 z_2 ( y_2 + z_2 )$ \\
  \hline
  7       & $y_2 z_5 = z_2 y_5$ \\
  \hline
  8       & $y_2 y_6 = y_3 y_5, z_2 y_6 = y_3 z_5$  \\
  \hline
  10      & $y_5^2 = y_2^2 y_6, z_5^2 = z_2^2 y_6$ \\
  \hline
  12      & $y_6^2 = y_2^2 z_2 y_6 + y_2 z_2^2 y_6 $ \\
  \hline
\end{tabular}  
\end{center}
We use the group structure to finish the calculations.  Fl\"oge's presentation for $\G5$ is $<A,B,M,R, S, U, W; A^2 = B^2 = M^2 = R^3 = S^3 = (AB)^2 = (BM)^2=1,\ AS=MR,\ U^{-1}AU = M, U^{-1}SU=R,\ W^{-1}MBW=AB,\ W^{-1}RW=S>$ \cite{Fl}.  We write this as a double HNN extension, with base group $G$, generated by $<A,B,M,R,S>$.  Consider the subgroups
\begin{eqnarray*}
  G_1 & = & <S,A,B; S^3=A^2=B^2=(AB)^2=1 > \ \cong \Z/3 \ast \Dt   \\
  G_2 & = & <R,B,M; R^3=B^2=M^2=(MB)^2=1 > \ \cong \Z/3 \ast \Dt   \\
  H   & = & <B,AS=MR; B^2 =1 >             \ \cong \Z/2 \ast {\Z}. \\
\end{eqnarray*}
Thus $G \cong G_1 {\ast}_H G_2 \cong (\Z/3 \ast \Dt) {\ast}_{(\Z/2 \ast {\Z}) } (\Dt \ast \Z/3)$.  Adding the group element $U$, set $G_3 = G \ast_{\Ps}$.  Adding $W$, $\Gamma_5 = G_3 \ast_{\Ps}$.  Both associated subgroups, $<A,S; A^2 = S^3 = 1>$ and $<MB,R; (MB)^2 = R^3 =1>$, are isomorphic to $\Ps$.

Use Formula \ref{LES} for an amalgamated product.
$$ \ldots \drarrow H^*(G;\Z)_{(2)} \brarrow H^*(\Dt;\Z)_{(2)} \oplus H^*(\Dt;\Z)_{(2)} \ararrow H^*(\Z/2 \ast \Z;\Z)_{(2)}  \drarrow \ldots.$$
By Lemma \ref{DtZtDt}, we see that $H^*(G;\Z)_{(2)} \cong H^*(\Dt \ast_{\Z/2} \Dt) \oplus \Z[\sigma_2]$; $\sigma_2$ is the image under $\delta$ of the torsion-free class in $H^1(\Z/2 \ast \Z;\Z)_{(2)}$.  Explicitly, $H^*(G;\Z)_{(2)} \cong \Z[x_2, y_2, z_2, y_3, z_3] \oplus \F[\sigma_2]$ with appropriate relations.  Assume that $y_2$ and $z_2$ are associated to $<A>$ and $<M>$ respectively.   

Next we add $U$ (twisting $A$ to $M$) to get $G_3$.  By Formula \ref{LES},    
$$ \ldots \drarrow H^*(G_3;\Z)_{(2)} \brarrow H^*(G;\Z)_{(2)} \ararrow H^*(\Ps;\Z)_{(2)}  \drarrow \ldots.$$
Denote the HNN class by $\sigma_1$.  The map $\alpha$ is surjective in positive even dimensions, and as $res^{\G3}_{\Ps} (y_2) = \theta^*(z_2)$, $y_2 + z_2$ is in the kernel of $\alpha$.  $H^2(G_3;\Z)_{(2)}$ contains two classes: $l_2=y_2+z_2$, and $m_2$, which is the preimage of $x_2$.  Also, $\alpha$ is the zero map in degree 3, so $H^3(G_3;\Z)_{(2)}$ contains the two classes corresponding to those in $H^3(G;\Z)_{(2)}$.  Finally, $H^*(G_3;\Z)_{(2)}$ has a new class in degree 4 (as $\alpha (x_2 y_2) = 0$) which we denote by $z_4$. 

The final calculation builds $\Gamma_5$ from $G_3$ by adding the group element $W$.  This element twists a copy of $\Ps$, sending $<MB>$ to $<AB>$.  Denote by $w_2$ the generator of $H^2(\Ps;\Z)_{(2)}$.  In terms of cohomology, $res^{\G5}_{\Ps}(x_2) = res^{\G5}_{\Ps}(z_2) = w_2$, and $\theta^* (x_2) = \theta^* (y_2) = w_2$.  Therefore, $\alpha$ sends both $l_2$ and $m_2$ to zero, so $\delta (w_2)$ has a non-trivial image in $H^3(\Gamma_5;\Z)_{(2)}$.  In addition, $\alpha$ is the zero map in degree 3.  Thus, there is a short exact sequence
$$ 0 \ararrow \Z/2 \drarrow H^3(\Gamma_5;\Z)_{(2)} \brarrow \Z/2 \oplus \Z/2 \ararrow 0.$$
The calculations from the Bockstein spectral sequence imply that there is only one torsion class in degree 3 of order 2.  We conclude that $H^3(\Gamma_5)$ consists of a class of order 2, and another class of order 4, which we call $t_3$.  We also get an HNN class, $\tau_1$.  That $t_3$ has trivial products follows from similar reasoning as in the calculation of $H^*(\Af \ast_{\Z/2} \Af ; \Z)_{(2)}$ in  Lemma \ref{AfZtAf}.  That $\tau_1$ has trivial products follows as it is the only class in the image of $\delta$.
To finish the argument, $\alpha$ is surjective in all other degrees.  For example, it is straightforward to check that $\alpha (z_4) \neq 0$.  Thus, the final long exact sequence breaks into short exact sequences, which implies that all other 2-torsion classes have order 2.  These have already been accounted for by the Bockstein spectral sequence, so the calculation is finished.


$H^*(\G5;\Z)_{(3)}$ is calculated in an identical manner as the case $\G6$.  \fin

\begin{theorem} $H^*(\G7;\Z)_{(2)} \cong \Z[y_2, \sigma_1]$ with relations $\sigma_1^2 = 2y_2 = 0$.  $H^*(\G7;\Z)_{(3)} \cong \Z[x_3,x_4] \oplus \Z[\sigma_1]$ with relations $3x_3 = 3x_4 = x_3^2 = 0$.  
\end{theorem}

\prf From \cite{Be}, $H^*(\G7;\F) \cong \F[x_1,y_1]$ with the single relation $y_1^2 = 0$.  The Bockstein spectral sequence yields one torsion class in the integral cohomology, $x_2 \leftrightarrow x_1^2 = Sq^1 (x_1)$, which is polynomial.  The integral class is harder to detect. 

We use an HNN description due to Fine, $\G7 \cong \left( \St \ast_{\Z/2} \St \right) \ast_{\Ps}$ \cite{F}.  The twisting sends a $\Z/2$ in the first copy of $\St$ to the second, and a $\Z/3$ in the second copy of $\St$ back to the first.   Note that $H^*(\St \ast_{\Z/2} \St;\Z)_{(2)} \cong H^*(\Z/2; \Z)_{(2)}$.  Furthermore, at the level of cohomology, the two dimensional generators of each copy of $H^*(\St;\Z)_{(2)}$ are identified in the amalgamated product.  

We use Formula \ref{LES} for an HNN extension.  
$$ \ldots \rightarrow H^*(\G7;\Z)_{(2)} \brarrow H^*(\St \ast_{\Z/2} \St;\Z)_{(2)} \ararrow H^*(\Ps;\Z)_{(2)}  \rightarrow \ldots$$
The restriction and twisting maps coincide, so $\alpha = 0$.  We denote the generator of $H^2(\Ps;\Z)_{(2)}$ by $w_2$, and the generator in $H^2(\G7;\Z)_{(2)}$ by $y_2$, and note that $res^{\G7}_{\Ps} (y_2) = w_2$.  There is also the HNN class, $\sigma_1$.  Substituting, the long exact sequence breaks into short exact sequences:
$$0 \rightarrow H^*(\Z/2;\Z) \drarrow H^{*+1}(\Gamma_7;\Z)_{(2)} \brarrow H^{*+1}(\St \ast_{\Z/2} \St;Z)_{(2)} \rightarrow 0. $$
The classes in the odd dimensions are the products $y_2^n \sigma_1$.  This follows from Lemma \ref{compat}, as $\delta (w_2^n) = \delta(res^{\G7}_{\Ps}(y_2^n) \cup 1) = y_2^n \cup \delta (1) = y_2^n \sigma_1$.   

At the prime 3, $H^*(\St \ast_{\Z/2} \St;\Z)_{(3)} \cong H^*(\Z/3;\Z) \widetilde{\oplus} H^*(\Z/3;\Z)$.  The rest of the argument is similar to the case $\G1$, with the exception of the HNN class.  Notice that no class in $H^*(\G7;\Z)_{(3)}$  restricts to the generator of $H^2(\Ps;\Z)_{(3)}$, so all products with the HNN class are trivial.  \fin

\begin{theorem} $H^*(\G{10} ; \Z)_{(2)} \cong \Z[y_2, z_2, y_3, y_5, z_5, y_6] \oplus \Z[t_3] \oplus \Z[\sigma_1] \oplus \Z[\tau_1] \oplus \Z[\eta_1] \oplus \Z[\sigma_2] \oplus \Z[\tau_2]$.  These satisfy the relations given in the table for $\G5$.  All torsion classes have order 2 except for $t_3$, which has order 4.  $H^*(\G{10};\Z) \cong \Z[x_2, \tau_1] \oplus \Z[\sigma_1] \oplus \Z[\eta_1] \oplus \Z[\sigma_2] \oplus \Z[\tau_2]$ with relations $\sigma_2^2 = \tau_1^2 = \eta_1^2 = \sigma_2^2 = \tau_2^2 = 3 x_2 = 0$.
\end{theorem}

\prf  The statement about the torsion follows as the polynomial classes in the mod-2 cohomology are identical to the case $\G5$.  The only difference is that $H^*(\G{10};\F)$ contains three exterior classes in each of degrees 1 and 2 that have trivial products with all other classes \cite{Be}.  

Fl\"oge's presentation for $\G{10}$ is $<A,B,L,S,D,U,W;\ A^2=B^2=L^2=S^3=(AB)^2=(AL)^2=1>$  with other relations involving $D,\ U$, and $W$ that we give in their respective extensions \cite{Fl}.  The base group $G_0$ has the presentation 
$<A,B,L,S;\ A^2=B^2=L^2=S^3=(AB)^2=(AL)^2=1>$ which we break into pieces.
\begin{eqnarray*}
  G_{01} & = & <A,B; A^2 = (AB)^2 = B^2 = 1 > \ \cong \Dt,   \\
  G_{02} & = & <A,L; A^2 = (AL)^2 = L^2 = 1 > \ \cong \Dt,   \\
  G_{03} & = & <S>    \ \cong \Z/3,    \\
  H      & = & <A>    \ \cong \Z/2.    
\end{eqnarray*}

With this decomposition, $G_0 \cong (G_{01} {\ast}_H\ G_{02})\ {\ast}\ G_{03} \cong (\Dt{\ast}_{\Z/2} \Dt) \ast \Z/3$.  The first HNN extension, $G_1$, adds the element $D$.  Its presentation is $<G_0, D; D^{-1}ALSD =  S^{-1}AB>$, with $<ALS> \ \cong {\mathbb Z}$.  $G_2$, the second HNN extension, adds the group element $U$.  Explicitly, $G_2 = \ <G_1, U;
U^{-1}DABD^{-1}U = D^{-1}ALD, U^{-1}LDS^{-1}D^{-1}U = BD^{-1}S^{-1}D>$.  The associated subgroup $<DABD^{-1}, LDS^{-1}D^{-1}>$ is isomorphic to $\Z/2 \ast {\Z}$.  The final HNN extension is $\G{10} = \ <G_2,W; W^{-1}BW = U^{-1}LU,  
W^{-1}D^{-1}SDW = U^{-1}DSD^{-1}U>$, where the associated subgroup $<B, D^{-1}SD>$ is isomorphic to $\Ps$. 

At the prime 2, $H^*(G_0;\Z)_{(2)} \cong H^*(\Dt \ast_{\Z/2} \Dt;\Z)_{(2)}$.  For the first HNN extension we use Formula \ref{LES}.
$$\ldots \drarrow H^*(G_1;\Z)_{(2)} \brarrow H^*(G_0;\Z)_{(2)} \ararrow H^*(\Z) \drarrow \ldots .$$
This implies that $H^*(G_1;\Z)_{(2)} \cong H^*(G_0;\Z)_{(2)}$ with two exceptions: there are new torsion-free classes in degrees 1 and 2.  All products with these classes vanish by Lemmas \ref{compat} and \ref{PY}.

The next two HNN extensions proceed in a similar manner as the case $\G5$.  These calculations account for the class $t_3$, and each HNN extension generates an HNN class in degree 1.  The second HNN extension also generates a torsion-free class in degree two, the image of $\delta H^1(\Z \ast \Z/2;\Z)$.  Therefore, $H^*(\G{10};\Z)_{(2)}$ has the same torsion classes as in the case of $\G5$, as well as three torsion-free classes in degree 1 and two in degree 2.  

In fact, $H^*(\G{10};\Z)_{(3)}$ is also close to $H^*(\G5;\Z)_{(2)}$.  Note that $H^*(G_0;\Z)_{(3)} \cong H^*(\Z/3;\Z)$.  These are the only torsion classes in $H^*(G_1;\Z)_{(3)}$ and $H^*(G_0;\Z)_{(3)}$, as these HNN extensions twist copies of $\Z$ and $\Z \ast \Z/2$.  The HNN extensions also give rise to two torsion-free classes in each of the degrees 1 and 2, the classes $\sigma_1$, $\sigma_2$, $\tau_1$, and $\tau_2$.  

It remains to determine the effect of the final HNN extension, where a copy of $\Ps$ is twisted.  The original copy of $\Z/3$ was generated by the group element $S$.  Although it is not twisted by $W$, a conjugate of it is, which induces an isomorphism on cohomology.  Therefore, we see that both the restriction and the twisting maps are isomorphisms onto $H^*(\Z/3;\Z)$, which implies $\alpha = 0$.  The long exact sequence breaks into short exact sequences, as in the case of  $\G5$.  We get an HNN class, $\eta_1 = \delta(1)$, and this class multiplies non-trivially with powers of the generator in $H^2(\G{10};\Z)_{(3)}$ that arises from $\alpha$ being the zero map.   \fin

\begin{theorem} $H^*(\G{11};\Z)_{(2)} \cong H^*(\Af \ast_{\Z/2} \Af;\Z)_{(2)} \oplus \Z[\sigma_1].$  The class $\sigma_1$ is torsion-free.  $H^*(\G{11};\Z)_{(3)} \cong \Z[x_2, \sigma_1]$ with relations $\sigma_1^2 = 3x_2 = 0$. 
\end{theorem}

\prf The Bockstein spectral sequence generates almost the entire integral ring structure from the mod-2 cohomology of the group, which is $H^*(\Gamma_{11};\F) \cong \F[u_2,\ v_3,\ \bar v_3,\ w_3] \oplus \F[s_1] \oplus \F[s_2]$ with relations $u_2^3 + w_3^2 + v_3^2 + \bar v_3^2 + w_3 (v_3 + \bar v_3) = 0,\ v_3 \bar v_3 = 0$.  With the exception of the class $\sigma_1$, this is identical to $H^*(\Af \ast_{\Z/2} \Af;\F)$.  The only class not accounted for is $\sigma_1$.

The alternate method of calculation resolves this issue.  Fine shows in \cite{F} that $\Gamma_{11} \cong \left( \Af {\ast}_{\Z/3} \Af \right) \ast_{\Ps}$, where the twisting sends copies of $\Z/2$ and $\Z/3$ from one $\Af$ to the other.
Thus, there is a long exact sequence, 
$$ \ldots \rightarrow H^*(\Gamma_{11};\Z)_{(2)} \brarrow H^*(\Af \ast_{\Z/3} \Af;\Z)_{(2)} \ararrow H^*(\Ps;\Z)_{(2)}  \rightarrow \ldots ,$$
with $H^*(\Ps;\Z)_{(2)} \cong H^*(\Z/2;\Z)_{(2)}$.   From this sequence we see that $\sigma_1$ is a torsion-free HNN class and the only class to arise from the image of $\delta$, so all products with $\sigma_1$ vanish.   $H^*(\Af \ast_{\Z/3} \Af;\Z)_{(2)} \cong H^*(\Af;\Z)_{(2)} \oplus H^*(\Af;\Z)_{(2)}$.  The result follows. 

For the 3-primary cohomology, start with $H^*(\Af \ast_{\Z/3} \Af;\Z)_{(3)} \cong H^*(\Z/3;\Z)$.  The twisting sends $\Z/3$ in one copy of $\Af$ to the other copy of $\Af$; when one calculates the cohomology of the HNN extension, this implies that $\alpha$ is zero.  We now proceed as in the case $\G2$.
\fin


\begin{remark}
Since the ring of integers $\calO_d$ contains number theoretic information, one might hope that some of this information would be detected in the cohomology ring of the corresponding Bianchi group $\G{d}$.  We note, however, that the torsion classes in $H^*(\G{5};\Z)$ and $H^*(\G{10};\Z)$ are identical--the difference is in the torsion-free classes!   It is not clear at this point if there is a connection, but if there is, it is clear that many more examples will need to be calculated to find it.
\end{remark}
   
\bibliographystyle{amsplain}

\begin{thebibliography}{99}

\bibitem{A-M} A. Adem and J. Milgram, \textit{Cohomology of Finite Groups}, 
Springer Verlag, 1994.



\bibitem{Be} E. Berkove, {\it On the Mod-2 Cohomology of the Bianchi Groups}, Trans. of the AMS.,
{\bf 352}, no. 10, (2000), 4585-4602.


\bibitem{Browd} W. Browder {\it Torsion in H-spaces}, Ann. of Math, (2) {\bf 74}, 1961, 24-51.

\bibitem{Bro} K. Brown, \textit{Cohomology of Groups}, Springer Verlag, 1993.

\bibitem{F} B. Fine, \textit{Algebraic Theory of the {B}ianchi Groups}, Marcel Dekker, Inc.
(Monographs and textbooks in pure and applied mathematics; 129), 1989.

\bibitem{Fl} D. Fl\"oge, \textit{Zur Struktur der $PSL_2$ \"uber einigen imagin\"ar-quadratischen Zahlringen}, Math. Zeitschrift \textbf{183} (1983), 255--279.

\bibitem{G-S} F. Grunewald and J. Schwermer, \textit{Subgroups of Bianchi Groups and Arithmetic Quotients of Hyperbolic 3-Space}, Trans. of the AMS (1) \textbf{335} (Jan. 1993), 47--78.

\bibitem{McC} J. McCleary, \textit{A User's Guide to Spectral Sequences}, Second Edition, Cambridge University Press (Cambridge Studies in advance mathematics; 58), 2001.

\bibitem{M} E. Mendoza, \textit{Cohomology of $PGL_2$ Over Imaginary Quadratic Integers}, Bonner Math. Schriften, vol. 128, 1980.

\bibitem{P} K. Pearson, \textit{Integral Cohomology and Detection of $\omega$-Basic 2-Groups}, Math. of Computation, Vol. 65 \textbf{213}, (Jan. 1996), 291-306.

\bibitem{P-Y} J. Pakianathan and E. Yalcin \textit{On G-Equivariant Coverings}, Preprint.

\bibitem{S-V} J. Schwermer and K. Vogtmann, \textit{The Integral Homology of $SL_2$ and $PSL_2$ of Euclidean Imaginary Quadratic Integers}, Comment. Math. Helvetica \textbf{58}  (1983), 573--598.

\bibitem{S} J. P. Serre, \textit{Trees}, Springer-Verlag, 1980. 

\bibitem{Sp} E. Spanier, \textit{Algebraic Topology}, 1st reprinted 
version, Springer-Verlag, 1989.


\bibitem{V} K. Vogtmann, \textit{Rational Cohomology of Bianchi Groups},
Math. Ann. \textbf{272} (1985), 399--419.

\end{thebibliography}

\end{document}